\documentclass[12pt]{article}
\usepackage{txfonts}

\usepackage{amssymb}
\renewcommand{\baselinestretch}{1.7}

\newtheorem{prop}{Proposition}
\newtheorem{cor}{Corollary}

\newtheorem{thm}{Theorem}
\newtheorem{defi}{Definition}

\newcommand {\beq}{\begin{equation}}
\newcommand {\eeq}{\end{equation}}
\newcommand {\beqa}{\begin{eqnarray}}
\newcommand {\eeqa}{\end{eqnarray}}         
\newcommand {\beqas}{\begin{eqnarray*}}
\newcommand {\eeqas}{\end{eqnarray*}}
\newcommand {\bea}{\begin{array}}
\newcommand {\eea}{\end{array}}
\newcommand {\bds}{\begin{displaymath}}
\newcommand {\eds}{\end{displaymath}}

\newcommand {\nn}{\nonumber}

\newcommand{\no}{\noindent}

\newcommand {\bebb}{\begin{thebibliography}}
\newcommand {\eebb}{\end{thebibliography}}      
\newcommand {\bbit}{\bibitem}

\def\dl{\delta}

\def\ot{\otimes}

\def\journal#1&#2(#3){\unskip, \sl #1\ \bf #2 \rm(19#3) }
\def\andjournal#1&#2(#3){\sl #1~\bf #2 \rm (19#3) }

\begin{document}

\begin{flushright}
\end{flushright}

\baselineskip =17pt

\vskip 1cm

\begin{center}
{\Large\bf On the Vertex Operators of the Elliptic \\
Quantum Algebra $U_{q,p}(\widehat{sl_2})_{k}$}

\vspace{1cm}

{\normalsize\bf Wen-Jing Chang $^{a, b}$ and Xiang-Mao Ding
$^a$\footnote{Partially Supported by NNSFC($\sharp$ 10671196;\,10231050);
corresponding author: E-mail: xmding@amss.ac.cn}}
\vspace{1cm}

{\em $^a$ Institute of Applied Mathematics, \\
Academy of Mathematics and Systems Science, \\
Chinese Academy of Sciences, P.O.Box 2734, Beijing, 100190, P.R.
China.}
\\
{\em $^b$ Graduate School of Chinese Academy of Sciences, China.}

\end{center}


\vskip 1.0cm

\begin{abstract}

A realization of the elliptic quantum algebra
$U_{q,p}(\widehat{sl_2})$ for any given level $k$ is constructed in
terms of three free boson fields and their accompanying twisted
partners. It can be viewed as the elliptic deformation of Wakimoto
realization. Two screening currents are constructed; they commute or
anti-commute with $U_{q,p}(\widehat{sl_2})$ modulo total
q-differences. The free fields realization for two types vertex
operators nominated as the type $I$ and the type $II$ vertex
operators are presented. The twisted version of the two types vertex
operators are also obtained. They all play crucial roles in
calculating correlation functions.

\end{abstract}


\renewcommand{\baselinestretch}{1.7}

\section{Introduction}

Infinite-dimensional symmetries, such as Virasoro algebra
($W$-algebra, for more general) and affine Lie algebra play a
central role in the two-dimensional Conformal Field Theories (2D
CFT) \cite{DiFMS}. While for the non-conformal (off-critical)
integrable theory, their roles are taken over by the so called
quantum algebras. From the algebraic point of view, there are three
kinds of quantum algebras, according to their different exchange
properties, which are nominated as rational, trigonometric and
elliptic quantum algebras respectively. The quantum algebras of the
former two kinds could be regarded as certain degenerate cases of
the last one. For example, the quantum affine algebra
(trigonometric), which is also known as the quantum group, and the
Yangian double with central (rational) could be obtained as a
certain limited case for the elliptic quantum algebras. Various
versions of elliptic quantum algebras, also known as elliptic
quantum groups \cite{Felder,Fronsdal,EF} have been introduced,
through an attempt to understand elliptic face models of statistical
mechanics, and in its semiclassical limit, a CFT of
Wess-Zumino-Witten (WZW) models on tori. Their roles are similar to
the Kac-Moody algebras in WZW models. From Hopf algebra point of
view, the elliptic quantum groups are nothing but quantum affine
algebras equipped with a co-product different from the original one
by a certain kind of twisting, so they can be viewed as quasi-Hopf
algebras in the sense of Drinfeld \cite{Drinfeld}. There are two
types of elliptic quantum groups which correspond to different types
of integrable models: the vertex type $A_{q,p}(\widehat{sl_N})$ and
the face type $B_{q,\lambda}(\mathcal{G})$, where $\mathcal{G}$ is a
Kac-Moody algebra associated with a symmetrizable generalized Cartan
matrix \cite{Kac}. The former is closely related to vertex models,
for example, the XYZ model, or equivalently, the eight vertex model
in the principal regime \cite{Baxter}; while some face models, such
as the Andrew-Baxter-Forrester (ABF) models \cite{ABF} which are
'solid-on-solid' (SOS) face models, possess symmetries corresponding
to the face type elliptic algebras.

From mathematical, it is natural to study these algebraic objects'
structures and their representations. In physical applications,
their representations are also required. The standard scheme to
study integrable models in field theories or statistical mechanics
is solving the following basic problems: to diagonalize the given
Hamiltonian, and then to compute the correlation functions. Usually,
it is quite difficult to solve such problems directly. While it has
been indicated that the algebraic analysis method is an extremely
powerful tool in studying solvable lattice models, especially in
deriving the correlation functions. This method is based on the
infinite dimensional quantum group symmetry possessed by a model and
the representation theory of such symmetry. As a result, if one
expect to perform algebraic analysis of the above both types of
elliptic lattice models, he should firstly study the corresponding
elliptic quantum groups and their representations.

In practice, the Wakimoto realization, or the so called free field
method, which is an infinite dimensional extension of the Heisenberg
algebra, is a quite effective and useful approach to study complicated
algebraic structures and their representations. The well known example
is the realization of the affine Lie algebras \cite{Wak,FF,Yu}. It is
also a common method to obtain representations of quantum affine algebras
\cite{Jimbo} and Yangian double \cite{Smir}. In \cite{DFJMN, JMMN},
the XXZ model in the anti-ferromagnetic regime was solved by
applying the level one representation theory of the quantum affine
algebra $U_q(\widehat{sl_2})$. In studying higher spin extension of
the XXZ model, free field realization at level $k>1$ is required,
which was constructed by several authors, such as \cite{Shiraishi,
Matsuo, DingW}. Furthermore, in \cite{AOS}, Wakimoto representation
of $U_q(\widehat{sl_N})$ with arbitrary level $k\geq 1$ was given,
and it plays a central role in understanding higher rank extension
of the XXZ model. Free field method is also a powerful way to study
the integrable massive field theory \cite{Luky}. Please see
\cite{JimMiw} for a nice review on quantum affine algebra, free
field realization and their applications. The level $k$ free field
representation of Yangian double $DY_{\hbar}(sl_2)$ and application
to physical problems were discussed in \cite{IK, Konno2}. The level
one free field realization of the Yangian double with central
$DY_{\hbar}(sl_N)$ was constructed in \cite{Iohara}; while the level
$k$ representation of $DY_{\hbar}(sl_N)$ and $DY_{\hbar}(gl_N)$ were
given in \cite{DHHZ}. It should also be remarked that the Yangian
double with central $DY_{\hbar}(\widehat{sl_2})$ is the symmetry
possessed by the Sine-Gordon model, which is the field theory limit
of the restricted SOS (RSOS) model \cite{DHZ,HY}. The bosonization
of the RSOS model was considered in \cite{Lukpagai}.

So following the algebraic approach, it is also important to obtain
the free fields realization of elliptic quantum algebras. For
example, in studying the RSOS model and its higher spin extension
(i.e. the $k$-fusion RSOS model), the free field representation of
$U_{q,p}(\widehat{sl_2})$ with any given level $k$ is needed, which
has been presented in \cite{Konno1} and this construction
corresponds to a deformation of special coset WZW model. In fact,
the elliptic algebra $U_{q,p}(\widehat{sl_2})$ is the Drinfeld
realization of $B_{q,\lambda}(\widehat{sl_2})$ showed in
\cite{JKOS}.

In fact, at classical level, there are various models of
representation for the current algebra, each one has its
significance in certain application. Here we just mention
two of them: free boson representation (Wakimoto construction)
\cite{Wak,FF}, and parafermion realization
\cite{Nemesch,GepQiu,Gep}. The similar realizations have been
extended to quantum affine algebras and Yangian double with central.
But there is no similar extension for the elliptic quantum algebras,
except the parafermion realization for $U_{q,p}(\widehat{sl_2})$ of
higher level $k$ \cite{Konno1}, and the parafermion representation
of $U_{q,p}(\widehat{sl_N})$ with level one \cite{KK}. In fact, the
realization presented in \cite{Konno1} is obtained by twisting the
parafermionic realization of the quantum affine algebra
$U_{q}(\widehat{sl_2})$, which can be considered as the elliptic
deformation of the parafermion realization. The $su(2)$ parafermion
currents can be identified with the coset WZW model of
$\widehat{sl(2)_k}\ot \widehat{sl(2)_1} / \widehat{sl(2)_{k+1}}$.
Although parafermion theory is important in physics
\cite{Nemesch,GepQiu,Gep} and in mathematics \cite{Dong}, it seems
that it cannot be used directly to study the elliptic quantum
algebra. In fact, the bosonization of non-local currents for higher
rank algebras is a huge project even in the classical level. So if
one want to deal with the elliptic quantum algebra of higher rank,
through bosonization of the non-local currents is not a practical
way.

It has special interest for the algebra of the intertwining
operators in the WZW model. It is derived by Knizhnik and
Zamolodchikov that the matrix coefficients of the intertwining
operators for the WZW model satisfy certain holonomic differential
equations, i.e., the Knizhnik-Zamolodchikov (KZ) equation \cite{KZ}.
There is an analogue holonomic q-difference equation for the quantum
affine vertex operators. They satisfy the quantum (q-deformed)
Knizhnik-Zamolodchikov (q-KZ) equation \cite{FR}. So it is also
expected that the representations of the elliptic quantum algebras
are likely to be helpful to construct the elliptic type solutions of
quantum Knizhnik-Zamolodchikov-Bernard (q-KZB) equation, which is a
higher genus extension of the q-KZ equation \cite{Bern}.
Furthermore, to consider a higher rank extension of the RSOS model,
we should construct free boson realization of
$U_{q,p}(\widehat{sl_N})$. However only in the level-one case,
parafermion realization of it was given in \cite{KK}, free boson
realization of it with higher level is unknown at present.

In this paper, we present a new free boson representation of
$U_{q,p}(\widehat{sl_2})$ with arbitrary level $k$. It is different
from the known one in \cite{Konno1} which was constructed in terms
of non-local currents. Our construction could be viewed as a twisted
version of the quantum semi-infinite flag manifolds \cite{FF}, which
is called the elliptic version of Wakimoto realization. The
realization of the quantum intertwining operators, such as the
screening currents and the vertex operators are also given. They are
necessary ingredients for calculating correlation functions and
investigating the irreducible representations. The screening
currents commute or anti-commute with $U_{q,p}(\widehat{sl_2})$, and
the integrations of such currents give the screening charges. For
$U_{q,p}(\widehat{sl_2})$, there are two types of intertwining
operators, which are called the type $I$ vertex operators (VO) and
the type $II$ vertex operators respectively with their different
physical significance. The former is a local operator which
describes the operation of adding one lattice site, and the formula
of the correlation functions can be expressed as traces of the
product of these operators over irreducible representation space;
while the latter play the role of particle creation or annihilation
operators. In fact, in this paper we construct two screening
currents, and the two types of VO's as well as their twisted ones.
Moreover we hope this construction can be generalized to other
cases, which will be considered in future study \cite{CD}.

The paper is organized as follows. In section $2$ we fix notations
and recall the definition of $U_{q,p}(\widehat{sl_2})$. We give our
construction of the free boson realization in section $3$. In
section $4$, two screening currents of it are constructed in terms
of free bosons. Finally in section $5$ we give the free boson
realization of the type $I$ VO's, the type $II$ VO's and their
twisted ones.

\setcounter{section}{1}
\setcounter{equation}{0}
\section{The definition of elliptic algebra $U_{q,p}(\widehat{sl_2})$}

Elliptic quantum algebras are introduced to study integrable models
with elliptic Boltzmann weights. There are two types of them: the
vertex type and the face type. Here we restrict to the face type
$B_{q,\lambda}(\mathcal{G})$ with $\mathcal{G}=\widehat{sl_2}$. For
the face type elliptic quantum algebras $U_{q,p}(\widehat{sl_2})$,
it can be considered as the Drinfeld realization of
$B_{q,\lambda}(\widehat{sl_2})$. In this section, we give a short
review on its definition.

Let us introduce a pair of parameters $p$ and $p^*$:
$$ p=q^{2r}, \ \
p^*=q^{2r^*}=pq^{-2c} \ \ \ \ \ \ \ (r^*=r-c; \ \
r,r^*\in\mathbb{R}_{>0})$$

\no here $c$ is the central element of the elliptic algebra
$U_{q,p}(\widehat{sl_2})$ defined below. Throughout this paper, the
complex number $q\neq0$, $|q|<1$ is fixed.

\no {\defi. The associative algebra $U_{q,p}(\widehat{sl_2})$ is
generated by the central element $c$ and the operator-valued
currents $H^{\pm}(z)$, $E(z)$ and $F(z)$ of the complex variable $z$
satisfying the following commutation relations:

\beqa
&&H^{\pm}(z)H^{\pm}(w)=(\frac{z}{w})^{2(\frac{1}{r^*}-\frac{1}{r})}
\frac{\Theta_p(q^{-2}\frac{z}{w})}
{\Theta_p(q^2\frac{z}{w})}\frac{\Theta_{p^*}(q^2\frac{z}{w})}
{\Theta_{p^*}(q^{-2}\frac{z}{w})}H^{\pm}(w)H^{\pm}(z), \\
&&H^{+}(z)H^{-}(w)=q^{2c(\frac{1}{r^*}+\frac{1}{r})}
(\frac{z}{w})^{2(\frac{1}{r^*}-\frac{1}{r})}
\frac{\Theta_p(pq^{-2-c}\frac{z}{w})}
{\Theta_p(pq^{2-c}\frac{z}{w})}\frac{\Theta_{p^*}(p^*q^{2+c}
\frac{z}{w})}{\Theta_{p^*}(p^*q^{-2+c}\frac{z}{w})}H^{-}(w)H^{+}(z), \\
&&H^{\pm}(z)E(w)=q^{\pm
\frac{c}{r^*}-2}(\frac{z}{w})^{\frac{2}{r^*}}
\frac{\Theta_{p^*}(q^{2\pm
\frac{c}{2}}\frac{z}{w})}{\Theta_{p^*}(q^{-2\pm\frac{c}{2}}
\frac{z}{w})}E(w)H^{\pm}(z), \\
&&H^{\pm}(z)F(w)=q^{\pm\frac{c}{r}+2}(\frac{z}{w})^{-\frac{2}{r}}
\frac{\Theta_p(q^{-2\mp\frac{c}{2}}
\frac{z}{w})}{\Theta_p(q^{2\mp\frac{c}{2}}\frac{z}{w})}F(w)H^{\pm}(z),
\\
&&[E(z), F(w)]=\frac{1}{(q-q^{-1})zw}[\dl(q^{-c}
\frac{z}{w})H^{+}(q^{-\frac{c}{2}}z)
-\dl(q^{c}\frac{z}{w})H^{-}(q^{-\frac{c}{2}}w)],
\\
&&E(z)E(w)=q^{-2}(\frac{z}{w})^{\frac{2}{r^*}}
\frac{\Theta_{p^*}(q^2\frac{z}{w})}
{\Theta_{p^*}(q^{-2}\frac{z}{w})}E(w)E(z), \\
&&F(z)F(w)=q^2(\frac{z}{w})^{-\frac{2}{r}}
\frac{\Theta_p(q^{-2}\frac{z}{w})}
{\Theta_p(q^2\frac{z}{w})}F(w)F(z),
 \eeqa}

\no where \beqas &&\dl(x)=\sum_{n\in\mathbb{Z}}x^n, \\
 &&\Theta_p(z)=(z;
p)_{\infty}(pz^{-1}; p)_{\infty}(p; p)_{\infty},
\\
&&(z; t_1,\cdots, t_k)_{\infty}=\prod_{n_1, \cdots,
n_k\geq0}(1-zt_1^{n_1}\cdots t_k^{n_k}), \eeqas

\no by the definition, $\Theta_p(z)$ is the standard elliptic
theta-function, up to a constant. In general, we can define
$\Theta_t(z)$ for any parameter $t=q^{2\nu}\ \ (\nu\in \mathbb{C})$
as
$$\Theta_t(z)=(z;
t)_{\infty}(tz^{-1}; t)_{\infty}(t; t)_{\infty}.$$ It should also be
remarked that this elliptic algebra $U_{q,p}(\widehat{sl_2})$
degenerates to the quantum affine algebra $U_q(\widehat{sl_2})$  in
the $p\rightarrow 0$ (or $r\rightarrow \infty$) limit.

In order to rewrite the relations (2.1)-(2.7) in more elegant form,
the following parameterization will be used: \beqas
&&q=e^{-\pi i/r\tau}, \\
&&p=e^{-2\pi i/\tau}, \ \ \ \ p^*=e^{-2\pi i/\tau^{*}}\\
&&z=q^{2u}=e^{-2\pi iu/r\tau}. \eeqas

\no In fact, we can use the notation of Jacobi theta function
$$\theta_{\nu}(u)=q^{\frac{u^2}{\nu}-u}
\frac{\Theta_{q^{2\nu}{}}(q^{2u})}{(q^{2\nu}; q^{2\nu})_{\infty}^3};$$
\no however for simplicity, we denote $\theta_{r}(u)$ as $\theta(u)$
and $\theta_{r^*}(u)$ as $\theta^*(u)$, which satisfy the
quasi-periodicity properties
$$\theta(u+r)=-\theta(u), \ \ \ \ \ \ \theta(u+r\tau)=-e^{-\pi \tau
i-2\pi iu/r}\theta(u),$$

\no and similar relations hold for $\theta^*(u)$ with $r$ replaced
by $r^*$. Then it is obvious to see that (2.1)-(2.7) can be
rewritten as follows:

\beqa &&H^{\pm}(u)H^{\pm}(v)=\frac{\theta(u-v-1)}
{\theta(u-v+1)}\frac{\theta^*(u-v+1)}{\theta^*(u-v-1)}H^{\pm}(v)H^{\pm}(u),
\\
&&H^{+}(u)H^{-}(v)=\frac{\theta(u-v-c/2-1)} {\theta(u-v-c/2+1)}
\frac{\theta^*(u-v+c/2+1)}
{\theta^*(u-v+c/2-1)}H^{-}(v)H^{+}(u),\\
&&H^{\pm}(u)E(v)=\frac{\theta^*(u-v\pm c/4+1)} {\theta^*(u-v\pm
c/4-1)}E(v)H^{\pm}(u),
\\
&&H^{\pm}(u)F(v)=\frac{\theta(u-v\mp c/4-1)}
{\theta(u-v\mp c/4+1)}F(v)H^{\pm}(u),\\
&&[E(u),
F(v)]=\frac{1}{(q-q^{-1})z w}[\dl(u-v-c/2)H^+(u-c/4)\nn\\
&& \hskip46mm
-\dl(u-v+c/2)H^-(v-c/4)],\\
&&E(u)E(v)=\frac{\theta^*(u-v+1)}{\theta^*(u-v-1)}E(v)E(u), \\
&&F(u)F(v)=\frac{\theta(u-v-1)}{\theta(u-v+1)}F(v)F(u), \eeqa

\no where the parameterizations $z=q^{2u}$ and $w=q^{2v}$ are
implicit in the above expressions. In the following, we will use
this parameterization without mentioning them if they do not make
confusion. Note that the above exchange relations have nice
periodicity property with the notations of the Jacobi theta
functions.

\setcounter{section}{2}
\setcounter{equation}{0}
\section{Fock Realization of $U_{q,p}(\widehat{sl_2})$ currents}

In this section, we construct an elliptic deformation of Wakimoto
realization for $U_{q,p}(\widehat{sl_2})$ using the free fields
representation of $U_q(\widehat{sl_2})$ for generic level $k$. We
first fix some conventions and review the Wakimoto realization of
$U_q(\widehat{sl_2})$; then give our construction of the realization
for $U_{q,p}(\widehat{sl_2})$ currents.

\subsection{Fock space and the quantum affine
algebra $U_{q}(\widehat{sl_2})$}

Three kinds free bosons a, b and c are needed to construct the
realization of $U_{q}(\widehat{sl_2})$. Their commutation relations
of modes are

\beqas &&[a_n,a_m]=\frac{[(k+2)n][2n]}{n}\dl_{n+m, 0},\ \ \ [p_a,
q_a]=2(k+2),\\
&&[b_n, b_m]=-\frac{[n]^2}{n}\dl_{n+m, 0},\ \ \ [p_b, q_b]=-1,\\
&&[c_n, c_m]=\frac{[n]^2}{n}\dl_{n+m, 0},\ \ \ [p_c, q_c]=1, \eeqas
and the others vanish, where $k$ is generic with $k\neq -2$.
Throughout this paper, the following standard symbol $[n]$ will be
used: $[n]=(q^n-q^{-n})/(q-q^{-1})$.

The vacuum state of the Fock space $|\mathbf{0}\rangle\equiv
|0,~0,~0\rangle$ is set as

$$a_{n}|\mathbf{0}\rangle=b_{n}|\mathbf{0}\rangle=c_{n}|\mathbf{0}\rangle=0
\ \ \ (n\geq0),$$

\no and a state $|l, m_1, m_2\rangle$ is produced through

$$|l, m_1, m_2\rangle\equiv
\exp\{l~q_a/2(k+2)+m_{1}q_b+m_{2}q_c\}|\mathbf{0}\rangle.$$

\no Obviously $|l, m_1, m_2\rangle$ is the highest weight state of
the bosonic Fock space, which is uniquely characterized by:

\beqas &&a_{n}|l, m_1, m_2\rangle=b_{n}|l, m_1, m_2\rangle=c_{n}|l,
m_1, m_2\rangle=0\
\ \ (n>0),\\
&&p_{a}|l, m_1, m_2\rangle=l|l, m_1, m_2\rangle,\ \ p_{b}|l, m_1,
m_2\rangle=-m_{1}|l, m_1, m_2\rangle,\\
&& p_{c}|l, m_1, m_2\rangle=m_{2}|l, m_1, m_2\rangle, \eeqas

\no then the Fock space ${\cal F}_{l, m_1, m_2}$ is generated by
negative modes $a_n$, $b_n$ and $c_n\ \ (n<0)$ acting on the highest
weight state $|l, m_1, m_2\rangle$. The dual Fock space could be
constructed  with the same matter.

For convenience, we denote free boson fields $a(z; \alpha)$  with
$\alpha \in\mathbb{C}$ and $a_{\pm}(z)$ as follows:

\beqas &&a(z; \alpha)=-\sum_{n\neq0}
\frac{a_n}{[n]}q^{-\alpha|n|}z^{-n}+q_a+p_a\ln z,\\
&&a_{\pm}(z)=\pm((q-q^{-1})\sum_{n>0}a_{\pm n}z^{\mp n}+p_a\ln q),
\eeqas

\no and $a(z; 0) \equiv a(z)$ for simplicity. Similarly, the free
boson fields $b(z; \alpha)$, $b_{\pm}(z)$ and $c(z; \alpha)$,
$c_{\pm}(z)$ can also be given. Normal order prescription $: \ :$ is
set by moving $a_n (n>0)$ and $p_a$ to right, while moving $a_n
(n<0)$ and $q_a$ to left. For example,
$$: \exp(a(z)) :=\exp(-\sum_{n<0}\frac{a_n}{[n]}z^{-n})
e^{q_a}z^{p_a} \exp(-\sum_{n>0}\frac{a_n}{[n]}z^{-n}).$$

With the help of the above free bosonic fields, four fields
$\psi_{\pm}(z)$ and $e^{\pm}(z)$ are introduced through their
actions on the Fock space ${\cal F}_{l, m_1, m_2}$. Let us fix the
actions of these currents on the Fock space as: $\psi_{\pm}(z)$:~
${\cal F}_{l, m_1, m_2}\mapsto {\cal F}_{l, m_1, m_2}$, $e^{+}(z):
{\cal F}_{l, m_1, m_2}\mapsto {\cal F}_{l, m_1-1, m_2-1}$ and
$e^{-}(z): {\cal F}_{l, m_1, m_2}\mapsto {\cal F}_{l, m_1+1,
m_2+1}$, respectively. Then these currents can be expressed as
follows:

\beqas &&\psi_{+}(z)=:\exp[b_{+}(q^{\frac{k}{2}}z)+a_{+}(qz)
+b_{+}(q^{\frac{k}{2}+2}z)]:,\\
&&\psi_{-}(z)=:\exp[b_{-}(q^{-\frac{k}{2}}z)+a_{-}(q^{-1}z)
+b_{-}(q^{-(\frac{k}{2}+2)}z)]:,\\
&&e^{+}(z)=\frac{-1}{(q-q^{-1})z}:\{\exp[b_+(z)-(b+c)(qz)]
-\exp[b_{-}(z)-(b+c)(q^{-1}z)]\}:,\\
&&e^{-}(z)=\frac{-1}{(q-q^{-1})z}:\{\exp[(b+c)(q^{-(k+1)}z)]
\exp[a_{-}(q^{-(\frac{k+2}{2})}z)+b_{-}(q^{-(k+2)}z)]\nn \\
&&\ \ \ \ \ \ \ \ \ \ \ \ \ \ \ \ \ \ \ \ \ \ \ \ \ \ \ \ \ \ \ \ \
-\exp[(b+c)(q^{k+1}z)]
\exp[a_{+}(q^{\frac{k+2}{2}}z)+b_{+}(q^{k+2}z)]\}:; \eeqas

\no and the following proposition is straightforward:

{\prop. The fields given above satisfy the following commutation
relations \cite{AOS}: \beqa
&&[\psi_{\pm}(z), \psi_{\pm}(w)]=0,\\
&&(z-q^{2-k}w)(z-q^{-2+k}w)\psi_{+}(z)\psi_{-}(w)\nn \\
&&\ \ \ \ \ \ \ \ =(z-q^{2+k}w)(z-q^{-2-k}w)\psi_{-}(w)\psi_{+}(z),\\
&&(z-q^{\pm (2-\frac{k}{2})}w)\psi_{+}(z)e^{\pm}(w) =(q^{\pm
2}z-q^{\mp \frac{k}{2}}w)e^{\pm}(w)\psi_{+}(z),\\
&&(z-q^{\pm (2-\frac{k}{2})}w)e^{\pm}(z)\psi_{-}(w)=(q^{\pm
2}z-q^{\mp \frac{k}{2}}w)\psi_{-}(w)e^{\pm}(z),\\
&&[e^{+}(z), e^{-}(w)]=\frac{1}{(q-q^{-1})z w}
[\delta(q^{-k}\frac{z}{w}) \psi_{+}(q^{-\frac{k}{2}}z)
-\delta(q^k\frac{z}{w})\psi_{-}(q^{-\frac{k}{2}}w)],\\
&&(z-q^{\pm 2}w)e^{\pm}(z)e^{\pm}(w)=(q^{\pm
2}z-w)e^{\pm}(w)e^{\pm}(z). \eeqa }

\no As a result, the above bosonic expression of the fields
$\psi_{\pm}(z)$ and $e^{\pm}(z)$ give the Wakimoto realization of
the quantum affine algebra $U_{q}(\widehat{sl_2})$ for generic level
$k$.

\subsection{Free boson realization of $U_{q,p}(\widehat{sl_2})$}

In this subsection, we will present a new free fields realization of
$U_{q,p}(\widehat{sl_2})$ with given level $k$. The elliptic algebra
$U_{q,p}(\widehat{sl_2})$ can be realized as the tensor product of
the elliptic currents $\Psi^{\pm}(z)$, $e(z)$, $f(z)$ of
$U_q(\widehat{sl_2})$ and a Heisenberg algebra \cite{JKOS}. The
elliptic currents $\Psi^{\pm}(z)$, $e(z)$ and $f(z)$ of
$U_q(\widehat{sl_2})$ are the fields satisfying the following
elliptic commutation relations:

\beqa
&&\Psi^{\pm}(z)\Psi^{\pm}(w)=\frac{\Theta_p(q^{-2}\frac{z}{w})}
{\Theta_p(q^2\frac{z}{w})}\frac{\Theta_{p^*}(q^2\frac{z}{w})}
{\Theta_{p^*}(q^{-2}\frac{z}{w})}\Psi^{\pm}(w)\Psi^{\pm}(z),
\\
&&\Psi^{+}(z)\Psi^{-}(w)=\frac{\Theta_p(pq^{-2-c}\frac{z}{w})}
{\Theta_p(pq^{2-c}\frac{z}{w})}\frac{\Theta_{p^*}(p^*q^{2+c}
\frac{z}{w})}{\Theta_{p^*}(p^*q^{-2+c}\frac{z}{w})}\Psi^{-}(w)\Psi^{+}(z),
\\
&&\Psi^{\pm}(z)e(w)=q^{-2}\frac{\Theta_{p^*}(q^{2\pm
\frac{c}{2}}\frac{z}{w})}{\Theta_{p^*}(q^{-2\pm\frac{c}{2}}
\frac{z}{w})}e(w)\Psi^{\pm}(z),
\\
&&\Psi^{\pm}(z)f(w)=q^2\frac{\Theta_p(q^{-2\mp\frac{c}{2}}
\frac{z}{w})}{\Theta_p(q^{2\mp\frac{c}{2}}\frac{z}{w})}f(w)\Psi^{\pm}(z),
\\
&&[e(z), f(w)]=\frac{1}{(q-q^{-1})zw}[\dl(q^{-c}\frac{z}{w})\Psi
^{+}(q^{-\frac{c}{2}}z)-\dl(q^{c}\frac{z}{w})\Psi^{-}(q^{-\frac{c}{2}}w)],
\\
&&e(z)e(w)=q^{-2}\frac{\Theta_{p^*}(q^2\frac{z}{w})}
{\Theta_{p^*}(q^{-2}\frac{z}{w})}e(w)e(z),\\
&&f(z)f(w)=q^2\frac{\Theta_p(q^{-2}\frac{z}{w})}
{\Theta_p(q^2\frac{z}{w})}f(w)f(z). \eeqa

In \cite{JKOS} the elliptic currents $\Psi^{\pm}(z)$, $e(z)$ and
$f(z)$ of $U_q(\widehat{sl_2})$ are realized by twisted the
parafermion realization of quantum affine algebra
$U_q(\widehat{sl_2})$. In the following, we will give another
bosonic realization, which can be viewed as the twisted quantum
version of the realization on flag manifold given by B. Feigin and
E. Frenkel \cite{FF}. To get the bosonic representation of these elliptic
currents, besides the bosonic fields introduced in the last
subsection, we need some new ones, such as $a_{\pm}^{*}(z)$,

\beqas
&&a_{+}^{*}(z)=-\sum_{n>0}\frac{a_n}{[rn]}z^{-n},\\
&&a_{-}^{*}(z)=\sum_{n>0}\frac{a_{-n}}{[r^*n]}z^n. \eeqas

\no Here we name them as the twisted partners of the fields
$a_{\pm}(z)$ respectively; and $b_{\pm}^{*}(z)$ are introduced with
the same matter by replacing $a_{\pm n}$ with $b_{\pm n}$ in the
above expressions. In terms of them we introduce two twisting
currents $U^{\pm}(z; r, r^*)$ depending on parameters $r$ and $r^*$
as:

\beqas &&U^{+}(z; r, r^*)=\exp[a_{-}^*(q^{r^{*}+\frac{k}{2}-1}z)
+b_{-}^*(q^{r^{*}-1}(q+q^{-1})z)],\\
&&U^{-}(z; r,
r^*)=\exp[a_{+}^*(q^{-(r-\frac{k}{2}-1)}z)+b_{+}^*(q^{-(r-k-1)}(q+q^{-1})z)].
\eeqas

By twisting the fields $\psi_{\pm}(z)$ and $e^{\pm}(z)$ in
subsection 3.1 with $U^{\pm}(z; r, r^*)$, we obtain the fields
$\Psi^{\pm}(z)$, $e(z)$ and $f(z)$ as:

\beqa
&&\Psi^{+}(z)=U^{+}(q^{\frac{k}{2}}z; r,
r^*)\psi_{+}(z)U^{-}(q^{-\frac{k}{2}}z; r, r^*),\\
&&\Psi^{-}(z)=U^{+}(q^{-\frac{k}{2}}z; r,
r^*)\psi_{-}(z)U^{-}(q^{\frac{k}{2}}z; r, r^*),\\
&&e(z)=U^{+}(z; r, r^*)e^{+}(z),\\
&&f(z)=e^{-}(z)U^{-}(z; r, r^*).
\eeqa

\no Obviously, the actions of the fields $\Psi^{\pm}(z)$, $e(z)$ and
$f(z)$ on the Fock space ${\cal F}_{l, m_1, m_2}$ are the same as
$\psi_{\pm}(z)$, $e^{+}(z)$ and $e^{-}(z)$ respectively, and we have
the following proposition:

{\prop.  The fields $\Psi^{\pm}(z)$, $e(z)$ and $f(z)$ obtained
above with $k=c$ satisfy the given elliptic commutation relations
(3.7)-(3.13).

Proof:} A straightforward but length OPE calculation verifies this
proposition. Here we just list some useful formulas:

\beqas &&e^{A}e^{B}=e^{[A, B]}e^{B}e^{A}, \ \ \ if\ [A, B]\ \
commute\ with\ A \ and\ B;\\
&&\exp(-\sum_{n>0}\frac{x^{n}}{n})=1-x;\\
&&(1-x)^{-1}=\sum_{n\geq 0}x^{n}. \ \ \ \ \ \ \ \ \ \Box\eeqas

By this proposition, we state that the currents in (3.14)-(3.17)
with $k=c$ give a bosonization of the elliptic currents of
$U_q(\widehat{sl_2})$. From their actions on the Fock space, we know
that all the currents keep the ``spin-$l/2$" representation.
Furthermore, it should also be remarked that, in the $p\rightarrow
0$ limit, $\Psi^{+}(z)\rightarrow (\psi_{-}(q^{k}z))^{-1}$,
$\Psi^{-}(z)\rightarrow (\psi_{+}(q^{k}z))^{-1}$, $e(z)\rightarrow
q^{-(2p_b+p_a)}(\psi_{-}(q^{k/2}z))^{-1}e^{+}(z)$ and
$f(z)\rightarrow e^{-}(z)q^{2p_b+p_a}(\psi_{+}(q^{k/2}z))^{-1}$ give
a new free fields realization of $U_q(\widehat{sl_2})$. It is
different from the one given in subsection 3.1.

However, the exchange relations of the currents given by the above
boson fields do not have good periodicity property. In order to
touch that goal, i.e., to construct a free boson realization of the
elliptic quantum algebra $U_{q,p}(\widehat{sl_2})$, we need
introduce a Heisenberg algebra generated by $\hat{p}$ and $\hat{q}$
such that
$$[\hat{q}, \hat{p}]=1,$$ and they commute with $a$, $b$ and $c$.
With the help of them we have new fields

\beqa &&H^{\pm}(u)=\Psi^{\pm}(z)e^{2\hat{q}}(q^{\mp
\frac{k}{2}}z)^{\frac{2p_b+p_a}{r}}
(q^{\pm(r-\frac{k}{2})}z)^{\frac{\hat{p}-1}{r}-\frac{\hat{p}-1}{r^*}},\\
&&E(u)=e(z)e^{2\hat{q}}z^{-\frac{\hat{p}-1}{r^*}},\\
&&F(u)=f(z)z^{\frac{2p_b+p_a}{r}}z^{\frac{\hat{p}-1}{r}}. \eeqa

\no To see them more clearly, we turn to the Fock space structure.
From the above expressions, the Fock space structure of
$H^{\pm}(u)$, $E(u)$ and $F(u)$ could be given as tensor product of
Fock spaces of $\Psi^{\pm}(z)$, $e(z)$ and $f(z)$ respectively, with
certain ones generated by $\hat{q}$. The results of their actions
are $H^{\pm}(u): {\cal F}_{l, m_1, m_2}\otimes{\cal F}_{n}\mapsto
{\cal F}_{l, m_1, m_2}\otimes{\cal F}_{n+2}$, $E(u): {\cal F}_{l,
m_1, m_2}\otimes{\cal F}_{n}\mapsto {\cal F}_{l, m_1-1,
m_2-1}\otimes{\cal F}_{n+2}$ and $F(u): {\cal F}_{l, m_1,
m_2}\otimes{\cal F}_{n}\mapsto {\cal F}_{l, m_1+1,
m_2+1}\otimes{\cal F}_{n}$. Here ${\cal F}_{n}$ is a trivial Fock
space generated by $|n\rangle\equiv e^{n\hat{q}}|\mathbf{0}\rangle$.
Then by applying proposition 2 and direct calculation, we can verify
the following theorem:

{\thm. The fields in Eqns. (3.18)-(3.20) with $k=c$ satisfy the
commutation relations (2.8)-(2.14).}

{\cor. $H^{\pm}(u)$, $E(u)$ and $F(u)$ given above realize the
elliptic algebra $U_{q,p}(\widehat{sl_2})$ with given level $k=c$.}

\setcounter{section}{3}
\setcounter{equation}{0}
\section{Construction of the screening currents}

In 2D CFT, screening current is a primary field of the
energy-momentum tensor with conformal weight $1$, and its
integration gives the screening charge. It has the property that it
commutes with the currents modulo a total differential of certain
field. This property ensures that the screening charge may be
inserted in the correlators by changing their conformal charges
without affecting their conformal properties. In this section, using
the bosons $a$, $b$ and $c$, we construct two screening currents
$S_{I}(z)$ and $S_{II}(z)$, which are integral parts in the free
fields approach. The two currents in this section commute or
anti-commute with the currents modulo a total q-difference of some
fields, so they could be regarded as a quantum deformation of the
screening currents in 2D CFT.

Denote a sort of q-difference operator with a parameter $n\in
\mathbb{Z}_{>0}$ by

$${}_n\partial_{z}X(z)\equiv\frac{X(q^{n}z)-X(q^{-n}z)}{(q-q^{-1})z},$$

\no which is called a total q-difference of a function $X(z)$.
Moreover to eliminate the total q-difference, one can define the
Jackson integral as

$$\int_{0}^{s\infty}X(z)d_{p}z\equiv s(1-p)\sum_{n\in
\mathbb{Z}}X(sp^{n})p^{n},$$

\no for a scalar $s\in \mathbb{C} \backslash \{0\}$ and a complex
number $p$ such that $\mid p \mid<1$. So that,

$$\int_{0}^{s\infty}({}_n\partial_{z}X(z))d_{p}z=0,$$

\no if it is convergent and take $p=q^{2n}$. For simplicity, we
denote boson fields with parameters $L_i$ and $M_j$ $(i, j \in
\mathbb{N})$ as follows:

\beqas &&A_{+}(L_1, \cdots, L_s; M_1, \cdots, M_{s+1}|z;
\alpha)=\sum_{n>0}\frac{[L_{1}n]\cdots[L_{s}n]}{[M_{1}n]
\cdots[M_{s+1}n]}a_{n}(q^{\alpha}z)^{-n},\\
&&A_{-}(L_1, \cdots, L_s; M_1, \cdots, M_{s+1}|z;
\alpha)=\sum_{n>0}\frac{[L_{1}n]\cdots[L_{s}n]}{[M_{1}n]
\cdots[M_{s+1}n]}a_{-n}(q^{\alpha}z)^{n}, \eeqas

\no and further abbreviate the notations as: \beqas &&A_{\pm}(L_1,
\cdots, L_s; M_1, \cdots, M_{s+1}|z)=A_{\pm}(L_1,
\cdots, L_s; M_1, \cdots, M_{s+1}|z; 0),\\
&&A_{\pm}(M|z; \alpha)=A_{\pm}(L_1, \cdots, L_s; L_1, \cdots,
L_s,M|z; \alpha);\eeqas

\no similarly the fields $B_{\pm}(L_1, \cdots, L_s; M_1, \cdots,
M_{s+1}|z; \alpha)$ and $C_{\pm}(L_1, \cdots, L_s; M_1, \cdots,
M_{s+1}|z; \alpha)$ can also be given.

Using these fields we obtain the screening currents as:

\beqas
&&S_{I}(z)=:\exp\{c(z)+\frac{q_{a}}{2}+q_{b}-r^{*}\hat{q}\}:,\\
&&S_{II}(z)=\frac{-1}{(q-q^{-1})z}:\exp\{A_{+}(k+2|z;
\frac{k+2}{2}) -\frac{1}{k+2}(q_{a}+p_a\ln z) \\
&&\hskip 4.25cm +A_{-}(-(k+2)|z;
-\frac{k+2}{2})\}\\
&&\hskip 3.5cm
 \times\{\exp[-b_{-}(z)-(b+c)(qz)]-\exp[-b_{+}(z)-(b+c)(q^{-1}z)]\}:, \eeqas

\no since they have the following properties:

{\thm: $S_I(z)$,~$S_{II}(z)$~satisfy the following relations with
the currents $H^{\pm}(z)$, $E(z)$ and $F(z)$ given by (3.18)-(3.20):

\beqas
&&\ \ \ \ \ \ \ H^{\pm}(z)S_{I}(w)=S_{I}(w)H^{\pm}(z)=O(1),\\
&&\ \ \ \ \ \ \ E(z)S_{I}(w)=-S_{I}(w)E(z)=
{}_{1}\partial_{w}[\frac{1}{z-w}\tilde{s}_{1}(z)]+O(1),\\
&&\ \ \ \ \ \ \ F(z)S_{I}(w)=-S_{I}(w)F(z)=O(1),\\
&&\ \ \ \ \ \ \ S_{I}(z)S_{I}(w)=-S_{I}(w)S_{I}(z)=O(1);\\
&&\ \ \ \ \ \ \ H^{\pm}(z)S_{II}(w)=S_{II}(w)H^{\pm}(z)=O(1),\\
&&\ \ \ \ \ \ \ E(z)S_{II}(w)=S_{II}(w)E(z)=O(1),\\
&&\ \ \ \ \ \ \ F(z)S_{II}(w)=S_{II}(w)F(z)={}_{(k+2)}
\partial_{w}[\frac{1}{z-w}\tilde{s}_{2}(z)]+O(1),\\
&&\ \ \ \ \ \ \
S_{II}(z)S_{II}(w)=\frac{\theta_{k+2}(u-v+1)}{\theta_{k+2}(u-v-1)}S_{II}(w)S_{II}(z);\\
&&\ \ \ \ \ \ \ S_{II}(z)S_{I}(w)=-S_{I}(w)S_{II}(z)={}_{1}
\partial_{w}[\frac{1}{z-w}\tilde{s}_{3}(z)]+O(1),
\eeqas

\no where the symbol $O(1)$ means regularity and $\tilde{s}_{i}(z)\
\ (i=1, 2, 3)$ are given by:

\beqas &&\tilde{s}_{1}(z)=:\exp\{A_{-}(r^*|z;
r-\frac{k}{2}-1)+\frac{q_a}{2}+B_{-}(-(r-k-2); r^*, 1|z;-1)\\
&&\hskip 2cm -p_{b}\ln z+B_{+}(1|z; -1)
+(2-r^*)\hat{q}-\frac{\hat{p}-1}{r^*}\ln z\}:,\\
&&\tilde{s}_{2}(z)=:\exp\{A_{-}(-(k+2)|z;
\frac{k+2}{2})+A_{+}(r-k-2;
k+2, r|z; \frac{k+2}{2})\\
&&\hskip 2cm -\frac{1}{k+2}(q_{a}+p_a\ln z) -B_{+}(2; 1, r|z;
-(r-k-1))
+\frac{2p_{b}+p_{a}+\hat{p}-1}{r}\ln z\}:,\\
&&\tilde{s}_{3}(z)=:\exp\{A_{-}(-(k+2)|z;
-\frac{k+2}{2})+A_{+}(k+2|z;
\frac{k+2}{2}) \\
&&\hskip 2cm +\frac{1}{2(k+2)}(kq_{a}-2p_a\ln z)
-b(z;1)+q_{b}-r^{*}\hat{q}\}:. \eeqas

Proof}:~~Straightforward calculation. Here we only take the last
relation as an example, we denote

$$S_{II}(z)\equiv\frac{-1}{(q-q^{-1})z}[A(z)-B(z)]$$

\no where \beqas &&A(z)=:\exp\{A_{-}(-(k+2)|z;
-\frac{k+2}{2})+A_{+}(k+2|z;
\frac{k+2}{2})\\
&&\hskip 2cm -\frac{1}{k+2}(q_{a}+p_a\ln z)-b_{-}(z)-(b+c)(qz)\}:,\\
&&B(z)=:\exp\{A_{-}(-(k+2)|z; -\frac{k+2}{2})+A_{+}(k+2|z;
\frac{k+2}{2})\\
&&\hskip 2cm -\frac{1}{k+2}(q_{a}+p_a\ln z)
-b_{+}(z)-(b+c)(q^{-1}z)\}:, \eeqas

\no then
$$S_{II}(z)S_{I}(w)=\frac{-1}{(q-q^{-1})z}[A(z)S_{I}(w)-B(z)S_{I}(w)],$$

\no and since the following relations hold: \beqas
&&A(z)S_{I}(w)=\frac{1}{qz-w}:A(z)S_{I}(w):,\ \ \ \mid z\mid
> \mid w\mid;\\
&&S_{I}(w)A(z)=\frac{1}{w-qz}:A(z)S_{I}(w):,\ \ \ \mid w\mid
> \mid z\mid;\\
&&B(z)S_{I}(w)=\frac{1}{q^{-1}z-w}:B(z)S_{I}(w):,\ \ \ \mid z\mid
> \mid w\mid;\\
&&S_{I}(w)B(z)=\frac{1}{w-q^{-1}z}:B(z)S_{I}(w):,\ \ \ \mid w\mid
> \mid z\mid,\eeqas

\no we obtain the following relation on the analytic continuations:

$$S_{II}(z)S_{I}(w)=-S_{I}(w)S_{II}(z)=\frac{-1}{(q-q^{-1})z}
[\frac{1}{qz-w}:A(z)S_{I}(w):-\frac{1}{q^{-1}z-w}:B(z)S_{I}(w):];$$

\no moreover,
$$:A(z)S_{I}(qz):=:B(z)S_{I}(q^{-1}z):\equiv\tilde{s}_{3}(z);$$

\no then by the definition of the total $q$-difference given above,
we get
$$S_{II}(z)S_{I}(w)=-S_{I}(w)S_{II}(z)={}_{1}
\partial_{w}[\frac{1}{z-w}\tilde{s}_{3}(z)]+O(1).\ \ \ \Box$$

\no It is easy to see that the actions of screening currents on the
Fock space are $S_{I}(z): {\cal F}_{l, m_1, m_2}\otimes{\cal
F}_{n}\mapsto {\cal F}_{l+(k+2), m_1+1, m_2+1}\otimes{\cal
F}_{n-r^{*}}$ and $S_{II}(z):{\cal F}_{l, m_1, m_2}\otimes{\cal
F}_{n}\mapsto {\cal F}_{l-2, m_1-1, m_2-1}\otimes{\cal F}_{n}$,
respectively. Please note that the current $S_{II}(z)$ acts
trivially on ${\cal F}_n$, so it is also the screening current of
the elliptic currents $\Psi^{\pm}(z)$, $e(z)$ and $f(z)$. On the
other side, the current $S_I(z)$ is not screening operator of them,
even if we remove the term involving $\hat{q}$ in $S_I(z)$ by hand.
In fact, using the expressions of $S_I(z)$ and $S_{II}(z)$ given
above, we can calculate the cohomology of the algebra and study the
irreducibility of modules of it, which will be discussed separately
in the future.

\setcounter{section}{4}
\setcounter{equation}{0}
\section{Realization of the Vertex Operators}

In fact, in 2D CFT, besides the screening currents, the other
important object that one should discuss is the primary field. In
WZW model, the primary fields could be realized as the highest
weight representation of Kac-Moody algebra, which are commonly known
as vertex operators (VOs)or intertwiner operators. They play crucial role in calculating correlation functions. For quantum affine algebra, there are two
types of vertex operators \cite{JimMiw} or intertwiner operators, in which the type $I$ is a local operator and could be regarded as the quantum counterpart
of the primary field in 2D CFT. In this section, we'll present a new
realization of the two types vertex operators and their twisted ones, which are different from the ones given in \cite{Konno1}, as they base on
distinct free fields realization of $U_{q,p}(\widehat{sl_2})$. For
the definitions of the VO's and the properties of them, please see
\cite{JimMiw, JKOS} for more details.

\subsection{The type $I$ and the type $II$ Vertex Operators}

For $U_{q,p}(\widehat{sl_2})$, the type $I$ vertex operators and the type $II$ vertex operators are defined to be the operators:

\beqa &&\widehat{\Phi_{l}}(u):\widehat{\cal
F}\rightarrow\widehat{\cal F}\otimes V_{l,
v}\\
&&\widehat{\Psi_{l}^{*}}(u):V_{l, v}\otimes\widehat{\cal
F}\rightarrow\widehat{\cal F} \eeqa

\no acting on the total Fock space $\widehat{\cal F}$, where $V_{l,
v}$ is the spin $\frac{l}{2}$ representation generated by vectors
$v_{m}^{l}\ \ (m=0, \cdots, l)$. For convenience, we set the
components $\Phi_{l, m}(u)$ and $\Psi_{l, m}^{*}(u)\ \ (m=0, \cdots,
l)$ of the VO's as

\beqas &&\widehat{\Phi_{l}}(u-\frac{1}{2})=\sum_{m=0}^{l}\Phi_{l,
m}(u)\otimes v_{m}^{l},\\
&&\widehat{\Psi_{l}^{*}}(u-\frac{k+1}{2})(v_{m}^{l}\otimes
\cdot)=\Psi_{l, m}^{*}(u). \eeqas

The fundamental property of the vertex operators is that they satisfy the
intertwining relations. In fact, intertwining operators of the algebra
could be used to define the Vertex operators in some sense. Here we just
pay our attention to the intertwining relations for the highest
components $\Phi_{l, l}(v)$ and $\Psi_{l, l}^{*}(v)$:

\beqa
&&H^{\pm}(u)\Phi_{l,l}(v)=\frac{\theta(u-v+\frac{l}{2}\mp\frac{k}{4})}
{\theta(u-v-\frac{l}{2}\mp\frac{k}{4})}\Phi_{l, l}(v)H^{\pm}(u),\\
&&E(u)\Phi_{l, l}(v)=\Phi_{l, l}(v)E(u),\\
&&F(u)\Phi_{l, l}(v)=\frac{\theta(u-v+\frac{l}{2})}
{\theta(u-v-\frac{l}{2})}\Phi_{l, l}(v)F(u);\\
&&H^{\pm}(u)\Psi_{l,
l}^{*}(v)=\frac{\theta^{*}(u-v-\frac{l}{2}\pm\frac{k}{4})}
{\theta^{*}(u-v+\frac{l}{2}\pm\frac{k}{4})}\Psi_{l, l}^{*}(v)H^{\pm}(u),\\
&&F(u)\Psi_{l, l}^{*}(v)=\Psi_{l, l}^{*}(v)F(u),\\
&&E(u)\Psi_{l, l}^{*}(v)=\frac{\theta^{*}(u-v-\frac{l}{2})}
{\theta^{*}(u-v+\frac{l}{2})}\Psi_{l,l}^{*}(v)E(u). \eeqa

\no It should be noted that all the expressions of the fields in
this section are considered to be normal-ordered.

Let us write currents $V^{\pm}(w; r,r^*)$ as:

\beqas &&V^{+}(w; r,r^*)=\exp\{-A_{+}(l, r^*; 2, k, r|w;
\frac{k+2}{2})-B_{+}(l, r^*; 1, k, r|w; k+1)\},\\
&&V^{-}(w; r,r^*)=\exp\{A_{-}(-l, r; 2, k, r^*|w;
\frac{k-2}{2})+B_{-}(-l, r; 1, k, r^*|w; -1)\}. \eeqas

\no Using them and the parameterization given in the second section,
we give a new realization of the type $I$ and the type $II$ VO's in
the following theorem:

{\thm. If we express $\Phi_{l, l}(v)$ and $\Psi_{l, l}^{*}(v)$
through:
\beqas
&&\Phi_{l, l}(v)=V^{+}(w; r,r^*)\exp\{A_{-}(l; 2,
k+2|w; \frac{k+2}{2})+A_{+}(l; k, k+2|w; \frac{k+2}{2})+B_{+}(l; 1,
k|w;1)\}\\
&&\hskip 3.1cm \times\exp\{\frac{l~q_a}{2(k+2)}
-\frac{l}{2r}(2p_{b}+\hat{p})\ln w\},\\
&&\Psi_{l, l}^{*}(v)=V^{-}(w; r,r^*)\exp\{-A_{-}(-l, k+1; 1, k,
k+2|w;
-\frac{k+2}{2})+A_{+}(-l; 2, k+2|w; \frac{k+2}{2})\\
&&\hskip 4.15cm -B_{-}(-l,
k+1; 1, 1, k|w; -1)+B_{+}(-l; 1, 1|w)\\
&&\hskip 4.15cm -C_{-}(-l;
1, 1|w)+C_{+}(-l; 1, 1|w)\}\\
&&\hskip 3.1cm \times\exp\{\frac{l~q_a}{2(k+2)}
+l(q_{b}+q_{c})-l\hat{q}+ \frac{l}{2r^*}\hat{p}\ln w\}, \eeqas

\no then they satisfy the intertwining relations (5.3)-(5.8).

Proof}:~~We take the relation (5.3) as an example. In fact, from the
bosonic expression of $\Phi_{l, l}(v)$ and (3.18), the following
OPE's can be derived:

\beqas &&H^{+}(u)\Phi_{l,l}(v)=q^{l-\frac{kl}{2r}}z^{\frac{l}{r}}
\frac{(q^{-l+\frac{k}{2}}\frac{w}{z};
p)_\infty}{(q^{l+\frac{k}{2}}\frac{w}{z};
p)_\infty}:H^{+}(u)\Phi_{l, l}(v):,\\
&&\Phi_{l,l}(v)H^{+}(u)=w^{\frac{l}{r}}
\frac{(pq^{-l-\frac{k}{2}}\frac{z}{w};
p)_\infty}{(pq^{l-\frac{k}{2}}\frac{z}{w};
p)_\infty}:H^{+}(u)\Phi_{l, l}(v):. \eeqas

\no The others can be derived similarly.$\ \ \ \ \ \ \Box$

Actually the lower components $\Phi_{l, m}(v)$ and $\Psi_{l,
m}^{*}(v)\ \ (m=0, \cdots, l)$ can be completely determined by the
highest ones $\Phi_{l, l}(v)$ and $\Psi_{l, l}^{*}(v)$, since they
obey the following recursive relations:

\beqa &&\Phi_{l,
m-1}(v)=F^{+}(v-\frac{l}{2})\frac{\theta(\hat{p}+h+l-m)}{\theta(\hat{p}+h)}\Phi_{l,
m}(v)\ \ \ \ \ (m=0, 1, \cdots, l),\\
&&\Psi_{l, m-1}^{*}(v)=\Psi_{l,
m}^{*}(v)E^{+}(v-\frac{l+k}{2}-r^*)\frac{\theta^{*}(m)\theta^{*}(\hat{p}-l+m-2)}{\theta^{*}(l-m+1)\theta^{*}(\hat{p}-2)}\
\ \ \ (m=0, 1, \cdots, l), \eeqa

\no where $E^{+}(v)$, $F^{+}(v)$ are half currents defined by

\beqas
&&E^{+}(v)=\varrho^{*}\oint_{C^{*}}E(v^{'})\frac{\theta^{*}(v-v^{'}+k/2-\hat{p}+1)
\theta^{*}(1)}{\theta^{*}(v-v^{'}+k/2)\theta^{*}(\hat{p}-1)}\frac{dw^{'}}{2\pi
iw^{'}},\\
&&F^{+}(v)=\varrho\oint_{C}F(v^{'})\frac{\theta(v-v^{'}+\hat{p}+h-1)\theta(1)}
{\theta(v-v^{'})\theta(\hat{p}+h-1)}\frac{dw^{'}}{2\pi
iw^{'}},\eeqas

\no and $h$ is one of the Drinfeld generators of
$U_{q}(\widehat{sl_2})$.

\no Here the contours are

\beqas
&&C^*: |p^{*}q^{k}w| < |w^{'}| < |q^{k}w|,\\
&&C: |pw| < |w^{'}| < |w|, \eeqas

\no and the constants $\varrho$, $\varrho^*$ are chosen to satisfy

$$\varrho\varrho^{*}\theta^{*}(1)\frac{\xi(q^{-2}; p^*, q)}{\xi(q^{-2}; p, q)}
 =q-q^{-1},$$

\no where the function $\xi(z; p, q)$ is

$$\xi(z; p, q)=\frac{(q^{2}z; p, q^4)_{\infty}(pq^{2}z; p, q^4)_{\infty}}
{(q^{4}z; p, q^4)_{\infty}(pz; p, q^4)_{\infty}}.$$

\subsection{The twisted Vertex Operators}

In this subsection, we discuss another two vertex operators ${\widehat{\Phi}_{l}}^{t}(u)$ and $\widehat{\Psi}_{l}^{*t}(u)$ for
$U_{q,p}(\widehat{sl_2})$, which are called the twisted type I VO's
and the twisted type II VO's (or twisted intertwiners) respectively.
Their definitions are analogously to the non-twisted ones. It means
that ${\widehat{\Phi}_{l}}^{t}(u)$ and $\widehat{\Psi}_{l}^{*t}(u)$
are also the operators of the same type as (5.1)-(5.2) and they also
have the similar decompositions, with their components denoted as
$\Phi_{l, m}^{t}(u)$ and $\Psi_{l, m}^{*t}(u)\ \ (m=0, \cdots, l)$.
However, the crucial difference between them lies in that
${\widehat{\Phi}_{l}}^{t}(u)$ and $\widehat{\Psi}_{l}^{*t}(u)$
satisfy the twisted intertwining relations. Here we also only
consider them for the highest components $\Phi_{l, l}^{t}(v)$ and
$\Psi_{l, l}^{*t}(v)$:

\beqa &&H^{\pm}(u)\Phi_{l,
l}^{t}(v)=\frac{\theta(u-v+\frac{l}{2}\mp\frac{k}{4})}
{\theta(u-v-\frac{l}{2}\mp\frac{k}{4})}\Phi_{l, l}^{t}(v)H^{\pm}(u),\\
&&E(u)\Phi_{l, l}^{t}(v)+\Phi_{l, l}^{t}(v)E(u)=0,\\
&&F(u)\Phi_{l, l}^{t}(v)=-\frac{\theta(u-v+\frac{l}{2})}
{\theta(u-v-\frac{l}{2})}\Phi_{l, l}^{t}(v)F(u);\\
&&H^{\pm}(u)\Psi_{l,l}^{*t}(v)
=\frac{\theta^{*}(u-v-\frac{l}{2}\pm\frac{k}{4})}
{\theta^{*}(u-v+\frac{l}{2}\pm\frac{k}{4})}\Psi_{l, l}^{*t}(v)H^{\pm}(u),\\
&&F(u)\Psi_{l, l}^{*t}(v)+\Psi_{l, l}^{*t}(v)F(u)=0,\\
&&E(u)\Psi_{l,l}^{*t}(v)=-\frac{\theta^{*}(u-v-\frac{l}{2})}
{\theta^{*}(u-v+\frac{l}{2})}\Psi_{l,l}^{*t}(v)E(u). \eeqa

\no It is easy to see that for the Cartan parts $H^{\pm}(u)$, there
is no difference for the twisted intertwining relations (5.3),
(5.11) and the non-twisted ones (5.6), (5.14); while for $E(u)$ and
$F(u)$, the difference between them is just a minus sign. With the
notations introduced before, we have

\beqas &&\Phi_{l, l}^{t}(v)=V^{+}(w; r,r^*)\exp\{A_{-}(k-l+1; 1,
k+2|w;
-\frac{k+2}{2}) +B_{-}(k-l+1; 1, 1|w; -1)\\
&&\hskip 3.75cm  -A_{+}(k-l, k+1; 1, k, k+2|w;
\frac{k+2}{2})-B_{+}(k-l, k+1; 1, 1, k|w; 1)\\
&&\hskip 3.75cm +C_{-}(k-l; 1, 1|w)-C_{+}(k-l; 1,
1|w)\}\\
&&\hskip 3.3cm \times
\exp\{\frac{(k-l)q_a}{2(k+2)}+(k-l)(q_{b}+q_{c})+r^{*}\hat{q}
+(\frac{l-r}{r}p_{c}-\frac{l}{2r}\hat{p})\ln w\},\\
&&\Psi_{l, l}^{*t}(v)=V^{-}(w; r,r^*)\exp\{A_{-}(-(k-l); k, k+2|w;
-\frac{k+2}{2})+B_{-}(-(k-l); 1, k|w; -1)\\
&&\hskip 3.95cm +A_{+}(l+2; 2, k+2|w; \frac{k+2}{2})+B_{+}(1|w;
l+1)\}\\
&&\hskip 3.3cm \times
\exp\{\frac{(k-l)q_a}{2(k+2)}-(l+r^*)\hat{q}+(-p_{b}+\frac{l}{2r^*}\hat{p})\ln
w\}. \eeqas

\no Similarly to the non-twisted case, we can prove that $\Phi_{l,
l}^{t}(v)$ and $\Psi_{l, l}^{*t}(v)$ obey the twisted intertwining
relations (5.11)-(5.16).

{\thm. $\Phi_{l, l}^{t}(v)$ and $\Psi_{l, l}^{*t}(v)$ gotten above
give a new realization of the twisted type $I$ and the twisted type
$II$ VO's. While the lower components $\Phi_{l, m}^{t}(v)$ and
$\Psi_{l, m}^{*t}(v)$ are determined by the same relations
(5.9)-(5.10) in which $\Phi_{l, m}(v)$ and $\Psi_{l, m}^{*}(v)$ are
replaced by $\Phi_{l, m}^{t}(v)$ and $\Psi_{l, m}^{*t}(v)$
respectively.}

In fact, there are a few degrees of freedom on choosing the
zero-modes of the above non-twisted and twisted vertex operators.
Here we only use the simplest ones.

\section{Discussion}

In this paper, we construct a new free boson realization of
$U_{q,p}(\widehat{sl_2})_k$ by twisting the flag manifold
realization, which can be viewed as the elliptic deformation of
Wakimoto realization. With this approach, the two important objects
(screening currents and Intertwiner Operators or vertex operators)
are also discussed in details. They all play important roles in
calculating correlation functions. Of course the derivation of
the multi-point correlation functions is a quiet interesting problem,
but in view of its complexity and the length of the manuscript, we will
discuss it in the future. Furthermore, it is an interesting problem to
extend our results to the other types of Lie algebras.

\section{Acknowledgments}

One of the authors (Ding)thanks the Kavli Institute for Theoretical
Physics China (KITPC) at the Chinese Academy of Sciences and the
program of String Theory and Cosmology for hospitality as the
write-up was completed. He is financially supported partly by the
Natural Science Foundations of China through the grands No.10671196
and No.10231050. He is also supported partly by a key Fund of
Chinese Academy of Sciences.

The authors also thank the referee of our paper, whose suggestions make
the expression in this paper more accurate and clearer.

\bebb{99}

\bbit{DiFMS}

P. Di Francesco, P. Mathieu, D. S\'{e}n\'{e}chal, {\it Conformal
Field Theory}, Springer, 1997.

\bbit{Felder}

G. Felder, Elliptic quantum groups,  Proc. ICMP Pairs 1994,
Cambridge-Hong Kong: International Press (1995) 211.

\bbit{Fronsdal}

C. Fr{\o}nsdal, Generalization and exact deformation of quantum
groups,  Publ. RIMS, Kyoto Univ, {\bf 33} (1997) 91.

\bbit{EF}

B. Enriquez, G. Felder, Elliptic quantum groups $E_{\tau,
\eta}(sl_2)$ and quasi-Hopf algebras,  Comm. Math. Phys. {\bf 195}
(1998) 651.

\bbit{Drinfeld}

V.G. Drinfeld, Quasi-Hopf algebras,  Leningrad Math. J. {\bf 1}
(1990) 1419.

\bibitem{Kac}

V.G. Kac, {\it Infinite dimensional Lie algebras}, third ed.,
Cambridge University press, 1990.

\bbit{Baxter}

R.J. Baxter, Partition function of the eight-vertex lattice model,
 Ann. Physics {\bf 70} (1972) 193.

\bbit{ABF}

G.E. Andrew, R.J. Baxter, P.J. Forrester, Eight-vertex SOS model and
generalized Rogers-Ramanujan-type identities,  J. Stat. Phys. {\bf
35} (1984) 193.

\bibitem{Wak}

M. Wakimoto, Fock representations of the affine Lie algebra
$A_{1}^{(1)}$,  Comm. Math. Phys. {\bf 104} (1986) 605.

\bbit{FF}

B. Feigin, E. Frenkel, Affine Kac-Moody algebras and semi-infinite
flag manifolds,  Comm. Math. Phys. {\bf 128} (1990) 161.

\bbit{Yu}

J.L. Petersen, J. Rasmussen,  M. Yu, Free field realizations of 2D
current algebras, screening currents and primary fields,  Nucl.
Phys. B{\bf 502} (1997) 649.

\bbit{Jimbo}

M. Jimbo, A q-difference analogue of $U~g$ and the Yang-Baxter
equation, Lett. Math. Phys. {\bf 10} (1986) 62.

\bbit{Smir}

F.A. Smirnov, Dynamical symmetries of massive integrable models $I$,
Int. J. Mod. Phys. A {\bf 7} suppl.1B (1992) 813; Dynamical
symmetries of massive integrable models $II$, Int. J. Mod. Phys.
A{\bf 7} suppl. {\bf 1B} (1992) 839.

\bbit{DFJMN}

B. Davies, O. Foda, M. Jimbo, T. Miwa, A. Nakayashiki,
Diagonalization of the XXZ Hamiltonian by vertex operators,  Comm.
Math. Phys. {\bf 151} (1993) 89.

\bbit{JMMN}

M. Jimbo, K. Miki, T. Miwa, A. Nakayashiki, Correlation functions of
the XXZ model for $\triangle <{-1}$,  Phys. Lett. A{\bf 168} (1992)
256.

\bbit{Shiraishi}

J. Shiraishi, Free boson representation of $U_q(\widehat{sl_2})$,
 Phys, Lett. A{\bf 171} (1992) 243.

\bbit{Matsuo}

A. Matsuo, Free field representation of quantum affine algebra
$U_q(\widehat{sl_2})$,  Phys. Lett. B{\bf 308} (1993) 260.

\bbit{DingW}

X.M. Ding, P. Wang, Parafermion representation of the quantum affine
$U_q(\widehat{sl_2})$, Modern. Phys. Lett. A{\bf 11} (1996) 921.

\bbit{AOS}

H. Awata, S. Odake, J. Shiraishi,  Free boson realization of
$U_q(\widehat{sl_N})$,  Comm. Math. Phys. {\bf 162} (1994) 61.

\bbit{Luky}

S. Lukyanov, Free field representation for massive integrable
models, Comm. Math. Phys. {\bf 167} (1995) 183.

\bbit{JimMiw}

M. Jimbo, T. Miwa,  Algebraic analysis of solvable lattice models,
 CBMS Regional Conference Series in Mathematics vol. {\bf 85} (1994) AMS
.

\bbit{IK}

K. Iohara, M. Kohno, A central extension of
$DY_{\hbar}(\widehat{sl_2})$ and its vertex representation,  Lett.
Math. Phys. {\bf 37} (1996) 319.

\bbit{Konno2}

H. Konno, Free field representation of level-$k$ Yangian double
$DY_{\hbar}(\widehat{sl_2})_{k}$ and deformation of Wakimoto
modules, Lett. Math. Phys. {\bf 40} (1997) 321.

\bbit{Iohara}

K. Iohara, Bosonic representations of Yangian double $DY_{\hbar}(g)$
with $g=gl_n$, $sl_n$,  J. Phys. A{\bf 29} (1996) 4593.

\bbit{DHHZ}

X.M. Ding, B.Y. Hou, B.Yuan. Hou, L. Zhao, Free boson representation
of $DY_{\hbar}({gl_N})_{k}$ and $DY_{\hbar}({sl_N})_{k}$,  J. Math.
Phys. {\bf 39} (1998) 2273.

\bbit{DHZ}

X.M. Ding, B.Y. Hou, L. Zhao, $\hbar$-Yangian deformation of the
Miura map and Virasoro algebra,  Intern. Jour. Mod. Phys. A{\bf 13}
(1998) 1129.

\bbit{HY}

B.Y. Hou, W.L. Yang, A $\hbar$-deformed Virasoro algebra as hidden
symmetry of the restricted sine-Gordon model, Comm. Theor. Phys.
{\bf 31} (1999) 265.

\bbit{Lukpagai}

S. Lukyanov, Y. Pugai, Multi-point Local Height Probabilities in the
Integrable RSOS Model, Nucl.Phys. B{\bf 473} (1996) 631.

\bbit{Konno1}

H. Konno, An elliptic algebra $U_{q,p}(\widehat{sl_2})$ and the
fusion RSOS model,  Comm. Math. Phys. {\bf 195} (1998) 373.

\bbit{JKOS}

M. Jimbo, H. Konno, S. Odake, J. Shiraishi, Elliptic algebra
$U_{q,p}(\widehat{sl_2})$: Drinfeld currents and vertex operators,
 Comm. Math. Phys. {\bf 199} (1999) 605.

\bbit{KK}

T. Kojima, H. Konno, The elliptic algebra $U_{q,p}(\widehat{sl_N})$
and the Drinfeld realization of the elliptic quantum group
$B_{q,\lambda}(\widehat{sl_N})$, Comm. Math. Phys. {\bf 239} (2003)
405.

\bbit{Nemesch}

D. Nemeschansky,  Feigin-Fuchs representation of $\widehat{su(2)}_k$
Kac-Moody algebras, Phys. Lett. B{\bf 224} (1989) 121.

\bbit{GepQiu}

D. Gepner, Z. Qiu, Modular invariant partition functions for
parafermionic field theories, Nucl. Phys. B{\bf 285} (1987) 423.

\bbit{Gep}

D. Gepner, New conformal field theories associated with Lie algebras
and their partition functions, Nucl. Phys. B{\bf 290} (1987) 10.

\bbit{Dong}

C.Y. Dong, J. Lepowsky, {\it Generalized Vertex Algebras and
Relative Vertex Operators}, Prog. Math. 112 Birh\"{a}user, 1993.

\bbit{KZ}

V.G. Knizhnik, A.B. Zamolodchikov, Current algebra and Wess-Zumino
models in two dimensions, Nucl. Phys. B{\bf 247} (1984) 83.

\bbit{FR}

I.B. Frenkel, N.Y. Reshetikhin, Quantum affine algebras and
holomorphic difference equations, Comm. Math. Phys. {\bf 146} (1992)
1 .

\bbit{Bern}

D. Bernard, On the Wess-Zumino-Witten model on the torus, Nucl.
Phys. B{\bf 303} (1988) 77.

\bbit{CD}

W.J. Chang, X.M. Ding, in preparing.

\eebb
\end{document}